\documentclass[12pt,a4paper]{article}
\usepackage[margin=1in]{geometry}
\usepackage{amssymb}
\usepackage{amsmath}
\usepackage{abstract}
\usepackage[utf8]{inputenc}
\usepackage[hidelinks]{hyperref}

\usepackage{titlesec}
\titleformat{\section}[block]{\center\scshape\large}{\thesection.}{0.7em}{}
\titleformat{\subsection}[block]{\flushleft}{\thesubsection.}{0.7em}{}
\renewcommand\thesection{\Roman{section}}

\usepackage{amsthm}
\newtheorem{theorem}{Theorem}
\newtheorem{lemma}{Lemma}

\theoremstyle{definition}

\title{On a weighted sum over multiplicative functions and its applications to the GPY sieve}

\author{Zihao Liu}
\date{}

\begin{document}
	\maketitle
	\begin{abstract}
		In this paper, we investigate the asymptotics of a class of weighted sums over multiplicative functions and apply our results to deduce a stronger asymptotic form of Yitang Zhang's smoothened GPY sieve with coefficient expressions that are friendly to numerical computations.
	\end{abstract}
	\section{Introduction}
	In number theory, asymptotic behaviors of sums in the form of
	\begin{equation}
		\label{eqn-mgdef}
		M_g(x,m,q)=\sum_{\substack{n\le x\\(n,q)=1}}g(n)\left(\log\frac xn\right)^m
	\end{equation}
	are often studied in order to investigate problems concerning primes. During the 19th and early 20th century, mathematicians such as Dirichlet \cite{dirichlet_ueber_2012}, F. Mertens \cite{mertens_ueber_1874}, and E. Landau \cite{landau_handbuch_1974} investigated the behavior of $M_g(x,0,q)$ when $g(n)$ is chosen to be $\tau_2(n)/n$, $\mu^2(n)/n$, or $\mu^2(n)/\varphi(n)$, resulting in asymptotic formulas that are valid when $x\to\infty$ and $q$ is fixed. In the 1950s, Harold N. Shapiro \& Jack Warga \cite{shapiro_representation_1950} and Yuan Wang \cite{wang_representation_1956} developed asymptotic estimates for sums of the type $M_g(x,0,q)$ that uniformly hold for all $q$ where $g(n)$ is chosen to be some multiplicative function supported on squarefree integers:
	\begin{align*}
		\sum_{\substack{n\le x\\(n,q)=1}}{\mu^2(n)\over n}={\varphi(q)\over q}\prod_{p\nmid q}\left(1-{1\over p^2}\right)\log x+O(\log\log q),
	\end{align*}
	\begin{align*}
		\sum_{\substack{n\le x\\(n,q)=1}}{\mu^2(n)2^{\omega(n)}\over n}
		&=\frac12\prod_{p|q}{p\over p+2}\prod_p\left(1-\frac1p\right)^2\left(1+\frac2p\right)\log^2x \\
		&+O(\log x\log\log xq)+O\{(\log\log q)^2\},
	\end{align*}
	\begin{align*}
		\sum_{\substack{n\le x\\(n,q)=1}}{\mu^2(n)\over\varphi(n)}={\varphi(q)\over q}\log x+O(\log\log q).
	\end{align*}
	Although the asymptotic expansions of these formulas exhibit various similarities, the derivation of each of them relies heavily on the Dirichlet convolution properties of specific arithmetical functions. The first general study in this direction was done by Ankeny and Onishi \cite{ankeny_general_1964} in 1964. Specifically, they applied a Tauberian argument\footnote{Proofs without use of Tauberian theorems are found in Chapter 6 of \cite{halberstam_sieve_1974} and Appendix A of \cite{friedlander_opera_2010}.} to prove that the general asymptotic relation
	\begin{equation}
		\label{eqn-ao}
		M_g(x,0,q)=\mathfrak S(q){(\log x)^k\over k!}+O_q\{(\log x)^{k-1}\}
	\end{equation}
	holds for fixed $q$ when $g(p)$ is well approximated by $k/p$ for some $k\ge1$ and $\mathfrak S(q)$ denotes the convergent product
	\begin{equation}
		\label{eqn-sq}
		\mathfrak S(q)=\prod_{p\nmid q}(1+g(p))\left(1-\frac1p\right)^k\prod_{p|q}\left(1-\frac1p\right)^k.
	\end{equation}
	Applying partial summation to \eqref{eqn-ao}, we see that when $m\ge0$ and $q$ is fixed there is
	\begin{equation}
		\label{eqn-opera}
		M_g(x,m,q)=\mathfrak S(q){m!\over(k+m)!}(\log x)^{k+m}+O_q\{(\log x)^{k+m-1}\}.
	\end{equation}
	
	In this article, we propose and prove \autoref{th-main}, which extends the validity of \eqref{eqn-opera} to the range $m\ge\max(0,-k)$ and offers an effective error term valid for all $x$ and $q$. In addition, we introduce and prove \autoref{th-smooth}, a smoothend version of \autoref{th-main}. After that, we apply these theorems to deduce a stronger form of Yitang Zhang's GPY sieve, which will be discussed in detail in \autoref{sc-zhang}.
	\section{Statement of results}
	\begin{theorem}
		\label{th-main}
		Let $q\in\mathbb Z^+$, $k\in\mathbb Z$, $\theta>0$, and $g(n)$ be a multiplicative function supported on squarefree numbers such that
		\begin{equation}
			\label{eqn-gp}
			g(p)=\frac kp+O\left(1\over p^{1+\theta}\right)\quad\forall p\nmid q,p\to+\infty.
		\end{equation}
		Then for all integer $m\ge\max(0,-k)$, there is
		\begin{align*}
			M_g(x,m,q)=\mathfrak S(q){m!\over(k+m)!}(\log x)^{k+m}+O\{(\log x)^{k+m-1}(\log\log q)^M\},
		\end{align*}
		where $M=M(k,m)\ge0$ is effectively computable and the implied constant in the O-term depends at most on $k$ and $m$.
	\end{theorem}

	\begin{theorem}
		\label{th-smooth}
		Let $M_g(x,m,q,z)$ be defined by
		\begin{equation}
			\label{eqn-mgsdef}
			M_g(x,m,q,z)=\sum_{\substack{n\le x\\(n,q)=1\\p|n\Rightarrow p<z}}g(n)\left(\log\frac xn\right)^m.
		\end{equation}
		If $u$ lies in a fixed interval, then under the assumption of \autoref{th-main}, we have
		\begin{align*}
			M_g(x,m,q,z)
			&=f(u;k,m)\mathfrak S(q){m!\over(k+m)!}(\log x)^{k+m} \\
			&+O\{(\log x)^{k+m-1}(\log\log q)^M\},
		\end{align*}
		where $f(u;k,m)$ is the solution to the differential-difference equation
		\begin{equation}
			\quad u^{k+m+1}f'(u)=-k(u-1)^{k+m}f(u-1)
		\end{equation}
		subjected to the initial condition that $f(u;k,m)=1$ for $0<u\le1$.
	\end{theorem}
	Since the $m=0$ case has already been proven by Ankeny and Onishi \cite{ankeny_general_1964}, $m>0$ is assumed throughout the rest of this paper. It is possible to extend the validity of these results to situations when $k$ and $m$ are not integers using Hankel contour and Gamma function, but the current version is sufficient for our purpose.
	\section{Lemmas}
	\begin{lemma}
		\label{lm-logs}
		For any $c>0$ and $m\in\mathbb Z^+$, we have
		\begin{align*}
			{m!\over2\pi i}\int_{c-i\infty}^{c+i\infty}{e^{ts}\over s^{m+1}}\mathrm ds
			=
			\begin{cases}
				t^m & t\ge0 \\
				0 & t<0
			\end{cases}.
		\end{align*}
	\end{lemma}
	\begin{proof}
		Move the path of integration leftward when $t<0$ and rightward otherwise.
	\end{proof}
	\begin{lemma}
		\label{lm-loglog}
		Let $q$ be any positive integer and $0\le\delta<1$. Then we have
		\begin{align*}
			\sum_{p|q}{\log p\over p^{1-\delta}}\ll(\log q)^\delta\log\log q,\quad\sum_{p|q}{1\over p^{1-\delta}}\ll(\log q)^\delta\log\log\log q,
		\end{align*}
		where the $\ll$ constant is absolute.
	\end{lemma}
	\begin{proof}
		Let $1<u<q$. Then we have
		\begin{align*}
			\sum_{p|q}{\log p\over p^{1-\delta}}
			&\le u^\delta\sum_{p\le u}{\log p\over p}+\sum_{\substack{p|q\\p>u}}{\log p\over u^{1-\delta}} \\
			&\ll u^\delta\log u+u^{\delta-1}\log q.
		\end{align*}
		Setting $u=\log q$ completes the proof. Proof for the second relation is similar except that we use the bound $\sum_{p\le u}p^{-1}\ll\log\log u$.
	\end{proof}
	\begin{lemma}
		\label{lm-zeta}
		There exists some $c_0>0$ such that for all $|v|\ge4$ and $u\ge1-{c_0\over\log|v|}$ there is
		\begin{align*}
			|\log\zeta(u+iv)|<\log\log|v|+O(1).
		\end{align*}
	\end{lemma}
	\begin{proof}
		See Theorem 6.7 of \cite{montgomery_multiplicative_2007}.
	\end{proof}
	\begin{lemma}
		\label{lm-buchstab}
		Let $S_g(x,q,z)$ be defined as follows
		\begin{align*}
			S_g(x,q,z)=\sum_{\substack{n\le x\\p|n\Rightarrow p<z\\(n,q)=1}}g(n)F\left(\log\frac xn\right).
		\end{align*}
		Then, for $2\le z\le x$, we have
		\begin{align*}
			S_g(x,q,z)=S_g(x,q,x)-\sum_{\substack{z\le p<x\\p\nmid q}}g(p)S_g\left(\frac xp,q,p\right).
		\end{align*}
	\end{lemma}
	\begin{proof}
		This is a direct generalization of Buchstab's identity in sieve theory.
	\end{proof}
	\begin{lemma}
		\label{lm-weight}
		Let $a_n$ be some sequence of complex numbers such that
		\begin{align*}
			\sum_{n\le x}a_n={A(\log x)^k\over k!}+O\{(\log x)^k\}.
		\end{align*}
		for some $A\in\mathbb C$ and $k\in\mathbb N$. Then for any $G\in C^1[0,1]$ we have
		\begin{align*}
			\sum_{n\le x}a_nG\left(\log x/n\over\log x\right)={A(\log x)^k\over(k-1)!}\int_0^1(1-t)^{k-1}G(t)\mathrm dt+O\{(\log x)^{k-1}\},
		\end{align*}
		where the implied constant depends on $A$, $k$, and $G$.
	\end{lemma}
	\begin{proof}
		See Lemma 4 of \cite{goldston_small_2009}.
	\end{proof}
	\section{The Dirichlet series associated with $g(n)$}
	Let $G(s)$ denote the Dirichlet series associated with $g(n)$:
	\begin{align*}
		G(s)=\sum_{\substack{n\ge1\\(n,q)=1}}{g(n)\over n^s}.
	\end{align*}
	In this section, we extract analytic properties of $G(s)$ in order to perform contour integration. From the Euler product formula, we can factor $G(s)$ into two parts:
	\begin{equation}
		\label{eqn-gsf}
		G(s)=\underbrace{\prod_{p\nmid q}\left(1+{g(p)\over p^s}\right)\left(1-{1\over p^{s+1}}\right)^k\prod_{p|q}\left(1-{1\over p^{s+1}}\right)^k}_{H(s)}\zeta^k(s+1).
	\end{equation}
	By \autoref{lm-zeta}, we see that $\zeta^k(s+1)\ll(\log|t|)^{|k|}$ whenever $\sigma\ge c_0/\log|t|,s=\sigma+it$. Thus, we only need to explore the properties of $H(s)$ in the rest of this section.
	\begin{lemma}
		\label{lm-conv}
		There exists an absolute constant $\delta>0$ such that for all $\sigma\ge-\delta$, the infinite product
		\begin{align*}
			\prod_{p\nmid q}\left(1+{g(p)\over p^s}\right)\left(1-{1\over p^{s+1}}\right)^k
		\end{align*}
		is absolutely convergent and uniformly bounded with respect to $q$.
	\end{lemma}
	\begin{proof}
		Plugging \eqref{eqn-gp} into the product, we see that
		\begin{align*}
			\log\left(1+{g(p)\over p^s}\right)\left(1-{1\over p^{s+1}}\right)^k
			&={g(p)-kp^{-1}\over p^s}+O\left(1\over p^{2\sigma+2}\right) \\
			&\ll{1\over p^{\sigma+1+\theta}}+{1\over p^{2\sigma+2}}.
		\end{align*}
		This indicates that the lemma's condition is satisfied when $0<\delta<\min(\theta,\frac12)$.
	\end{proof}
	\begin{lemma}
		\label{lm-logq}
		There exists an absolute constant $C>0$ such that for any $\sigma\ge-\delta$, we have
		\begin{align*}
			|H(s)|\ll(\log\log q)^{C|k|(\log q)^\delta}.
		\end{align*}
	\end{lemma}
	\begin{proof}
		Due to \autoref{lm-conv}, it suffices to show that
		\begin{align*}
			\prod_{p|q}\left(1-{1\over p^{s+1}}\right)^k\ll(\log\log q)^{C|k|(\log q)^\delta}.
		\end{align*}
		Taking logarithm, we have
		\begin{align*}
			\left|\log\prod_{p|q}\left(1-{1\over p^{s+1}}\right)^k\right|\le|k|\sum_{p|q}{1\over p^{1-\delta}}+O(1),
		\end{align*}
		so plugging \autoref{lm-loglog} into the right hand side completes the proof.
	\end{proof}
	\section{Proof of \autoref{th-main}}
	Setting $t=\log x/n$ in \autoref{lm-logs} and plugging it into \eqref{eqn-mgdef}, we have for $c>0$ that
	\begin{align*}
		M_g(x,m,q)={m!\over2\pi i}\int_{c-i\infty}^{c+i\infty}\sum_{\substack{n\ge1\\(n,q)=1}}g(n)\left(\frac xn\right)^s{\mathrm ds\over s^{m+1}}={m!\over2\pi i}\int_{c-i\infty}^{c+i\infty}G(s){x^s\over s^{m+1}}\mathrm ds.
	\end{align*}
	To estimate the integral, we first truncate the path of integration so that it follows from \autoref{lm-logq} that for any large $T\le x$ there is
	\begin{equation}
		\label{eqn-mgi}
		M_g(x,m,q)={m!\over2\pi i}\int_{c-iT}^{c+iT}G(s){x^s\over s^{m+1}}\mathrm ds+O\left\{x^c(\log\log q)^{C|k|}c^{-|k|}\over T^m\right\}.
	\end{equation}
	In order to estimate $M_g(x,m,q)$ using residue theorem, we apply \autoref{lm-logq} so that when $\delta_0=c_1(\log T)^{-1}(\log\log q)^{-1}$ for some $c_1=\min(c_0,\log2)$ there is
	\begin{equation}
		\label{eqn-iv}
		\int_{-\delta_0-iT}^{-\delta_0+iT}G(s){x^s\over s^{m+1}}\mathrm ds\ll x^{-\delta_0}(\log\log q)^{2C|k|}.
	\end{equation}
	\begin{equation}
		\label{eqn-ih}
		\left(\int_{c-iT}^{-\delta_0-iT}+\int_{-\delta_0+iT}^{c+iT}\right)G(s){x^s\over s^{m+1}}\mathrm ds\ll x^c(\log\log q)^{2C|k|}{(\log T)^{|k|}\over T^{m+1}}.
	\end{equation}
	Now, setting $c=1/\log x$ and $\log T=\sqrt{\log x}$ converts the above bounds of \eqref{eqn-mgi}, \eqref{eqn-iv}, and \eqref{eqn-ih} into
	\begin{equation}
		\label{eqn-mg}
		M_g(x,m,q)
		={m!\over2\pi i}\oint_{(0+)}G(s){x^s\over s^{m+1}}\mathrm ds+O\{(\log\log q)^{2C|k|}\}.
	\end{equation}
	Therefore, the remaining task is to calculate the residue of $G(s)x^s/s^{m+1}$ at $s=0$. For convenience, let $I(s)$ be defined by
	\begin{equation}
		\label{eqn-idef}
		I(s)=s^kG(s)=H(s)[s\zeta(s+1)]^k.
	\end{equation}
	Then it follows from \eqref{eqn-gsf} that $I(s)$ is analytic near $s=0$ and from \eqref{eqn-sq} that $I(0)=H(0)=\mathfrak S(q)$. These facts allow us to transform the residue integral of \eqref{eqn-mg} into
	\begin{align}
		\label{eqn-mg-main}
		{1\over2\pi i}\oint_{(0+)}I(s){x^s\over s^{k+m+1}}\mathrm ds
		&=\mathfrak S(q){(\log x)^{k+m}\over(k+m)!} \\
		\label{eqn-mg-err}
		&+\sum_{1\le n\le k+m}{(\log x)^{k+m-n}\over(k+m-n)!}{1\over2\pi i}\oint_{(0+)}{I(s)\over s^{n+1}}\mathrm ds.
	\end{align}
	Evidently, \eqref{eqn-mg-main} will become the main term in the expansion of $M_g(x,m,q)$, so the remaining task is to estimate \eqref{eqn-mg-err}. By \autoref{lm-logq}, \eqref{eqn-idef}, and Cauchy's inequality, we observe that when $r=\log2/\log\log q$ there is
	\begin{align*}
		\oint_{|s|=r}{I(s)\over s^{n+1}}\mathrm ds\ll(\log\log q)^{2C|k|+n}.
	\end{align*}
	Plugging this back into \eqref{eqn-mg-err} we obtain
	\begin{equation}
		\label{mg-errb}
		\sum_{1\le n\le k+m}\ll(\log x)^{k+m-1}(\log\log q)^{2C|k|+k+m}.
	\end{equation}
	Finally, combining \eqref{mg-errb} and \eqref{eqn-mg} completes the proof of \autoref{th-main} with $M=2C|k|+k+m$.
	\section{Proof of \autoref{th-smooth}}
	The case where $0<u\le1$ follows directly from \autoref{th-main}, so it suffices to prove for $u\ge1$.

	Assume \autoref{th-smooth} is true when $0<u\le r-1$ for some integer $r\ge2$, so when $r-1<u\le r$ it follows from \autoref{lm-buchstab} that
	\begin{equation}
		\label{mg-buchstab}
		M_g(x,m,q,x^{1/u})=M_g(x,m,q)-\underbrace{\sum_{\substack{x^{1/u}\le p<x\\p\nmid q}}g(p)M_g\left(\frac xp,m,q,p\right)}_\Sigma.
	\end{equation}
	Because ${\log x/p\over\log p}\le u-1\le r-1$, we apply the inductive hypothesis to deduce
	\begin{align*}
		\Sigma
		&=\mathfrak S(q){m!\over(k+m)!}(\log x)^{k+m}\underbrace{\sum_{\substack{x^{1/u}\le p<x}}\frac kpf\left({\log x\over\log p}-1;k,m\right)\left(1-{\log p\over\log x}\right)^{k+m}}_\Gamma \\
		&+O\left\{\mathfrak S(q)(\log x)^{k+m}\sum_{x^{1/u}\le p<x}{1\over p^{1+\theta}}\right\}+O\left\{\mathfrak S(q)(\log x)^{k+m}\sum_{\substack{x^{1/u}\le p<x\\p|q}}\frac1p\right\} \\
		&+O\left\{(\log\log q)^M(\log x)^{k+m-1}\sum_{x^{1/u}\le p<x}\frac1p\right\}.
	\end{align*}
	By \autoref{lm-loglog} and some standard estimates, we have
	\begin{align*}
		\mathfrak S(q)\ll(\log\log q)^{|k|},\quad\sum_{x^{1/u}\le p<x}\frac1p\ll1,\quad\sum_{x^{1/u}\le p<x}{1\over p^{1+\theta}}\ll x^{-\theta/u},
	\end{align*}
	and
	\begin{align*}
		\sum_{\substack{x^{1/u}\le p<x\\p|q}}\frac1p\ll(\log x)^{-1}\sum_{p|q}{\log p\over p}\ll(\log x)^{-1}\log\log q.
	\end{align*}
	Moreover, by partial summation, we have
	\begin{align*}
		\Gamma
		&=k\int_{x^{1/u}}^x{\mathrm dt\over t\log t}f\left({\log x\over\log t}-1\right)\left(1-{\log t\over\log x}\right)^{k+m}+O\left(1\over\log x\right) \\
		&=k\int_1^uf(v-1;k,m)(1-v^{-1})^{k+m}{\mathrm dv\over v}+O\left(1\over\log x\right).
	\end{align*}
	Plugging all of these results back into \eqref{mg-buchstab}, we obtain
	\begin{align*}
		M_g(x,m,q,x^{1/u})
		&=\left\{1-k\int_1^uf(v-1;k,m)(v-1)^{k+m}v^{-k-m-1}\mathrm dv\right\} \\
		&\times\mathfrak S(q){m!\over(k+m)!}(\log x)^{k+m}+O\{(\log\log q)^M(\log x)^{k+m-1}\}.
	\end{align*}
	For $r-1<u\le r$, if we define $f(u;k,m)$ by
	\begin{align*}
		f(u;k,m)=1-k\int_1^uf(v-1;k,m)(v-1)^{k+m}v^{-k-m-1}\mathrm dv,
	\end{align*}
	the differentiating both side with respect to $u$ yields the differential-difference equation stated in \autoref{th-smooth}, thus completing the proof.
	\section{Yitang Zhang's smoothened GPY sieve}\label{sc-zhang}
	A positive real number $\theta$ is called the level of distribution of primes in arithmetic progressions if for every $\varepsilon>0$ and every $A>0$,
	\begin{equation}
		\label{eqn-eh}
		\sum_{d\le x^{\theta-\varepsilon}}\max_{(a,d)=1}\left|\pi(x;q,a)-{\pi(x)\over\varphi(d)}\right|\ll_A{x\over\log^Ax}\quad(x\ge2)
	\end{equation}
	Under Bombieri--Vinogradov theorem \cite{davenport_multiplicative_1980}, \eqref{eqn-eh} holds unconditionally for $\theta=\frac12$. For the sake of brevity, we use $\operatorname{EH}(\theta)$ to denote the hypothesis that \eqref{eqn-eh} holds for some given $\theta$. It is conjectured by Elliott and Halberstam \cite{elliott_conjecture_1970} that $\operatorname{EH}(1)$ holds, and Bombieri--Vinogradov theorem is equivalent to $\operatorname{EH}(\frac12)$.

	In 2005, Goldston, Y\i ld\i r\i m, and Pintz \cite{goldston_primes_2009} proved that if $\operatorname{EH}(\theta)$ holds for some $\theta>\frac12$, then there will be infinitely many pairs of primes such that their absolute difference is bounded by some constant that only depends on $\theta$. In other words, there exists some $C(\theta)\ge2$ such that
	\begin{equation}
		\label{eqn-gap}
		\liminf_{n\to\infty}(p_{n+1}-p_n)\le C(\theta)
	\end{equation}

	Let $\chi_\mathbb P$ denote the characteristic function for primes. Then the GPY sieve is expressed as follows
	\begin{equation}
		\label{eqn-gpy}
		S=\sum_{N<n\le2N}\left(\sum_{1\le i\le k}\chi_\mathbb P(n+h_i)-1\right)w_n^2.
	\end{equation}
	where $h_1<h_2<\dots<h_k$ is chosen such that $Q(n)=\prod_{1\le i\le k}(n+h_i)$ is not always divisible by a fixed prime divisor\footnote{These numbers are called an admissible $k$-tuple. A trivial example is to let $h_1,h_2,\dots,h_k$ be the first $k$ primes greater than $k$.}, and $w_n$ is some well chosen sequence of real numbers. If $S>0$ for some $N$, then there exists some $N<n\le2N$ and some $1\le r,s\le k$ such that $n+h_r$ and $n+h_s$ are both primes, so $p_{n+1}-p_n\le h_k-h_1$ for infinitely many $n$.

	In the original work of Goldston, Y\i ld\i r\i m, and Pintz, the shape of $w_n$ is as follows:
	\begin{equation}
		\label{eqn-gpy-wn}
		w_n=\sum_{d\in D_n}\lambda_d,\quad D_n=\{d\le R:d|Q(n)\},\quad \lambda_d=\mu(d)P\left(\log R/d\over\log R\right),
	\end{equation}
	In their original paper, the authors plugged in $P(x)=x^m$ for $m>k$ and obtained
	\begin{align*}
		S=[1+o(1)]C_1(k,m,\theta){N\over\log^kN}.
	\end{align*}
	In addition, they showed that if $\operatorname{EH}(\theta)$ holds for some $\theta>\frac12$, then there exists some $k$ and $m$ such that $C_1(k,m,\theta)>0$. In 2007, Soundararajan \cite{soundararajan_small_2007} showed that even when the choice of $P(x)$ is loosened to $x^{-k}P(x)\in C[0,1]$, it is still necessary to assume $\operatorname{EH}(\theta)$ for some $\theta>\frac12$ in order that $C_1(k,m,\theta)$ is positive. Therefore, it is sufficient and necessary to assume $\operatorname{EH}(\theta)$ for some $\theta>\frac12$ for the prototypical GPY sieve to produce bounded gaps between primes.

	In 2014, Yitang Zhang \cite{zhang_bounded_2014} introduced smoothing to the GPY sieve. While still using $P(x)=x^m$ for $\lambda_d$, he restricted the size of prime factors of integers in $D_n$:
	\begin{equation}
		\label{eqn-zhang}
		D_n=\{d\le R:d|Q(n)\wedge p|d\Rightarrow d<z\},\quad z=N^\delta
	\end{equation}
	This choice of $D_n$ corresponds to a smoothened variant of hypothesis $\operatorname{EH}(\theta)$. Let $\operatorname{EH}(\theta,\delta)$ denote the hypothesis that for every $A>0$,

	\begin{equation}
		\label{eqn-eh2}
		\sum_{\substack{d\le x^\theta\\ p|d\Rightarrow p<N^\delta}}\max_{(a,d)=1}\left|\pi(x;q,a)-{\pi(x)\over\varphi(d)}\right|\ll_A{x\over\log^Ax}\quad(x\ge2).
	\end{equation}

	Instead of looking for exact asymptotic, Zhang estimated the difference between smoothened GPY sieve and the prototypical GPY sieve to deduce a lower bound:
	\begin{equation}
		\label{eqn-zhang-bound}
		S\ge C_2(k,m,\theta,\delta){N\over\log^kN}.
	\end{equation}
	Zhang proved $\operatorname{EH}(\theta,\delta)$ when $\delta=1/1168$ and $\theta=1/2+2\delta$ and that under this choice of parameters, $C_2(k,m,\theta,\delta)>0$ for $k=3.5\times10^6$, $m=k+180$, and sufficiently large $N$, thereby deducing the existence of infinitely many pairs of primes that are $<7\times10^7$ apart.

	Zhang did not derive an asymptotic formula in \cite{zhang_bounded_2014} thus unable to utilize the full capacity of the smoothened GPY sieve. In this section, we apply \autoref{th-smooth} to deduce an asymptotic form of Zhang's original smoothened GPY sieve so that it would be possible for one to analyze its strengths and limits. Moreover, this would allow one to apply this sieve to number theory questions other than prime gaps. Apart from \autoref{th-smooth} itself, the derivation is purely elementary. The asymptotic coefficients obtained using this approach is also friendly to numerical computations as they are defined using iteration formulas.

	\begin{theorem}
		\label{th-zhang}
		Let $S$ be defined as in \eqref{eqn-gpy} such that
		\begin{align*}
			w_n=\sum_{\substack{d\le R\\d|Q(n)\\p|d\Rightarrow p\le z}}\lambda_d,\quad\lambda_d=\mu(d)P\left(\log R/d\over\log R\right),\quad P(x)=x^m,\quad m>k
		\end{align*}
		Then under $\operatorname{EH}(\theta,\delta)$, we have for $z=N^\delta$ and $R=N^{\theta/2}$ that
		\begin{align*}
			S\sim{N\over(\log R)^k}\prod_p\left(1-{\nu_p\over p}\right)\left(1-\frac1p\right)^{-k}\left\{{k\theta\over2}I_{k-1}(1,u)-I_k(1,u)\right\},\quad u={\theta\over2\delta}
		\end{align*}
		where $\nu_d$ is the number of roots of $Q(n)$ in $\mathbb Z/d\mathbb Z$ and
		\begin{align*}
			I_s(t,v)=\int_0^1{(1-x)^{s-1}\over(s-1)!}[f(uxt;-s,m)P^{(s)}(xt)]^2\mathrm dx,\quad0<v\le1
		\end{align*}
		\begin{align*}
			I_s(t,v)=I_s(t,1)-s\int_1^vI_s[(1-x^{-1})t,x-1](1-x^{-1})^s{\mathrm dx\over x}.\quad v>1
		\end{align*}
		in which $P^{(s)}(x)$ denotes the $s$'th derivative of $P(x)$.
	\end{theorem}

	We deduce \autoref{th-zhang} from the following auxiliary result:

	\begin{theorem}
		\label{th-sieve}
		Let $k$ be some positive integer and $g(n)$ be some multiplicative function supported on squarefree numbers such that
		\begin{align*}
			g(p)=\frac kp+O\left(1\over p^2\right).
		\end{align*}
		Then we have for $z=R^{1/u}$ that
		\begin{align*}
			Q_g=\sum_{\substack{d_1,d_2\le R\\p|[d_1,d_2]\Rightarrow p<z}}\lambda_{d_1}\lambda_{d_2}g([d_1,d_2])\sim{I_k(1,u)\over(\log R)^k}\prod_p(1-g(p))\left(1-\frac1p\right)^{-k}.
		\end{align*}
	\end{theorem}
	\renewcommand*{\proofname}{Proof of \autoref{th-zhang}}
	\begin{proof}
		By \eqref{eqn-gpy}, we see that $S=\sum_{1\le i\le k}S_i-S_0$, where
		\begin{align*}
			S_0=\sum_{N<n\le2N}w_n^2,\quad S_i=\sum_{N<n\le2N}\chi_\mathbb P(n+h_i)w_n^2.
		\end{align*}
		Thus, it suffices to evaluate $S_0$ and $S_i$ separately. By a change of order of summation, we have
		\begin{align*}
			S_0=\sum_{\substack{d_1,d_2\le R\\p|[d_1,d_2]\Rightarrow p<z}}\lambda_{d_1}\lambda_{d_2}\sum_{\substack{N<n\le 2N\\Q(n)\equiv0([d_1,d_2])}}1=NQ_{g_0}+E_0.
		\end{align*}
		where $g_0(d)=\nu_d/d$, which is multiplicative due to Chinese remainder theorem and
		\begin{align*}
			E_0
			&\ll\sum_{\substack{d\le R^2\\p|d\Rightarrow p<z}}\mu^2(d_1)\mu^2(d_2)\nu_{[d_1,d_2]}\le R^2\sum_{\substack{d\\p|d\ge1\Rightarrow p<z}}\mu^2(d)\tau_3(d){\nu_d\over d} \\
			&=R^2\prod_{p<z}\left(1+{3\nu_p\over p}\right)\le R^2\prod_{p<z}\left(1+{3k\over p}\right)\ll R^2(\log z)^{3k}.
		\end{align*}
		
		On the other hand, it follows from the prime number theorem that
		\begin{align*}
			S_i
			&=\sum_{\substack{d_1,d_2\le R\\p|[d_1,d_2]\Rightarrow p<z}}\lambda_{d_1}\lambda_{d_2}\sum_{\substack{N<n\le2N\\Q(n)\equiv0([d_1,d_2])}}\chi_\mathbb P(n+h_i) \\
			&=Q_{g_1}\sum_{N<n\le2N}\chi_\mathbb P(n+h_i)+E_i \\
			&=[1+o(1)]{N\over\log N}Q_{g_1}+E_i,
		\end{align*}
		where $g_1$ is multiplicative such that $g_1(p)=(\nu_p-1)/\varphi(p)$ and
		\begin{align*}
			E_i
			&\ll\sum_{\substack{d\le R^2\\p|d\Rightarrow p<z}}\mu^2(d)\tau_{3k-3}(d)\max_{(a,d)=1}\left|\sum_{\substack{N<t\le2N\\t\equiv a(d)}}\chi_\mathbb P(n+h_i)-{1\over\varphi(d)}\sum_{N<t\le2N}\chi_\mathbb P(n+h_i)\right|,
		\end{align*}
		which is $\ll_A N/\log^AN$ by an application of Cauchy--Schwarz inequality and $\operatorname{EH}(\theta,\delta)$.

		By the properties of polynomials, we know that $Q(n)$ has exactly $k$ roots in $\mathbb Z/p\mathbb Z$ if and only if $p$ does not divide its discriminant, a finite number. Therefore, when $p\to+\infty$ we have
		\begin{equation}
			\label{eqn-g0g1}
			\nu_p=k,\quad g_0(p)=\frac kp,\quad g_1(p)={k-1\over p}+O\left(1\over p^2\right).
		\end{equation}
		Finally, plugging \eqref{eqn-g0g1} into \autoref{th-sieve}, so we have
		\begin{equation}
			\label{eqn-qg0}
			S_0\sim N\cdot{I_k(1,u)\over(\log R)^k}\prod_p\left(1-{\nu_p\over p}\right)\left(1-\frac1p\right)^{-k}
		\end{equation}
		and for $1\le i\le k$
		\begin{equation}
			\label{eqn-qg1}
			S_i\sim{N\over\log N}\cdot{I_{k-1}(1,u)\over(\log R)^{k-1}}\prod_p\left(1-{\nu_p-1\over p-1}\right)\left(1-\frac1p\right)^{-(k-1)}.
		\end{equation}
		Notice that
		\begin{align*}
			\left(1-{\nu_p-1\over p-1}\right)=(p-\nu_p)(p-1)^{-1}=\left(1-{\nu_p\over p}\right)\left(1-\frac1p\right)^{-1},
		\end{align*}
		so combining \eqref{eqn-qg0} and \eqref{eqn-qg1} proves the asymptotic formula for $S$.
	\end{proof}
	\renewcommand*{\proofname}{Proof}
	Although in this paper we only use \autoref{th-smooth} to refine Zhang's sieve, we expect that by modifying the methods applied to prove \autoref{th-sieve} to higher rank situation (i.e. Maynard--Tao sieve), thereby potentially improving current bounds for prime gaps.
	\section{Proof of \autoref{th-sieve}}
	Let $h(l)$ be a multiplicative function such that $h(p)=g(p)/(1-g(p))$. Then by M\"obius inversion, we have
	\begin{align*}
		Q_g
		&=\sum_{\substack{d_1,d_2\le R\\p|[d_1,d_2]\Rightarrow p<z}}\lambda_{d_1}g(d_1)\lambda_{d_2}g(d_2){1\over g((d_1,d_2))}=\sum_{\substack{d_1,d_2\le R\\p|[d_1,d_2]\Rightarrow p<z}}\lambda_{d_1}g(d_1)\lambda_{d_2}g(d_2)\sum_{l|(d_1,d_2)}{1\over h(l)} \\
		&=\sum_{\substack{l\le R\\p|l\Rightarrow p<z}}{\alpha_l^2\over h(l)},\quad\alpha_l=\sum_{\substack{d\le R\\d\equiv0(l)\\p|d\Rightarrow p<z}}g(d)\lambda_d=\mu(l)g(l)\sum_{\substack{t\le R/l\\(t,l)=1\\p|t\Rightarrow p<z}}\mu(t)g(t)\left(\log R/lt\over\log R\right)^m.
	\end{align*}
	For convenience, let $A$ denote the following quantity
	\begin{equation}
		\label{eqn-gp}
		A=\prod_p(1-g(p))\left(1-\frac1p\right)^k.
	\end{equation}
	Then it follows that when $q$ is squarefree,
	\begin{align*}
		\prod_{p\nmid q}(1+\mu(p)g(p))\left(1-\frac1p\right)^k\prod_{p|q}\left(1-\frac1p\right)^k=A\prod_{p|q}(1-g(p))^{-1}=A{h(q)\over g(q)}.
	\end{align*}
	Therefore, when we apply \autoref{th-smooth}, we have
	\begin{align*}
		\alpha_l
		&={\mu(l)g(l)\over(\log R)^m}M_{\mu g}\left(\frac Rl,m,l,z\right) \\
		&={A\mu(l)h(l)\over(\log R)^m}f\left({\log R/l\over\log z};-k,m\right){m!\over(m-k)!}\left(\log\frac Rl\right)^{m-k}+O\left\{g(l)(\log R/l)^{m-k-1}\log\log l\over(\log R)^m\right\} \\
		&={A\mu(l)h(l)\over(\log R)^k}f\left(u\cdot{\log R/l\over\log R};-k,m\right)P^{(k)}\left(\log R/l\over\log R\right)+O\left\{g(l)(\log R)^{m-k-1+\varepsilon}\over(\log R)^m\right\}.
	\end{align*}
	For convenience, let $G(x)=[f(ux;-k,m)P^{(k)}(x)]^2$, so we have
	\begin{equation}
		\label{eqn-alpha}
		{\alpha_l^2\over h(l)}={h(l)A^2\over(\log R)^{2k}}G\left(\log R/l\over\log R\right)+O\left\{g(l)\over(\log R)^{2k+1-\varepsilon}\right\}.
	\end{equation}
	Using the multiplicative property of $g$, we have
	\begin{align*}
		\sum_{\substack{l\le R\\p|l\Rightarrow p<z}}g(l)
		&\le\sum_{p|l\Rightarrow p<z}g(l)=\prod_{p<z}(1+g(p))\le\exp\left\{\sum_{p<z}\frac kp+O\left(\sum_p{1\over p^2}\right)\right\} \\
		&=\exp\{k\log\log z+O(1)\}\ll(\log z)^k\ll(\log R)^k.
	\end{align*}
	Therefore, plugging \eqref{eqn-alpha} into $Q_g$ gives
	\begin{equation}
		\label{eqn-qgh}
		Q_g={A^2\over(\log R)^{2k}}\sum_{\substack{l\le R\\p|l\Rightarrow l<z}}h(l)G\left(\log R/l\over\log R\right)+O\left\{1\over(\log R)^{k+1-\varepsilon}\right\}.
	\end{equation}
	For convenience, we let $H(x,z)$ denote the following quantity:
	\begin{equation}
		\label{eqn-hxz}
		H(x,z)=\sum_{\substack{l\le x\\p|l\Rightarrow l<z}}h(l)G\left(\log x/l\over\log R\right).
	\end{equation}
	Note that $h(p)=g(p)/(1-g(p))=k/p+O(1/p^2)$ and
	\begin{align*}
		\prod_p(1+h(p))\left(1-\frac1p\right)^k=\prod_p(1-g(p))^{-1}\left(1-\frac1p\right)^k=\frac1A,
	\end{align*}
	so it follows from \autoref{th-main} that
	\begin{align*}
		\sum_{l\le R}h(l)={A^{-1}\over k!}(\log R)^k+O\{(\log R)^{k-1}\}.
	\end{align*}
	Plugging this into \autoref{lm-weight}, we see that for $2\le x\le R$ there is
	\begin{align*}
		H(x,x)
		&=\sum_{l\le x}h(l)G\left({\log x/l\over\log x}\cdot{\log x\over\log R}\right) \\
		&=A^{-1}(\log x)^k\int_0^1{(1-t)^{k-1}\over(k-1)!}G\left(t\log x\over\log R\right)\mathrm dt+O\{(\log x)^{k-1}\}. \\
		&=A^{-1}(\log x)^kI_k\left({\log x\over\log R},1\right)+O\{(\log x)^{k-1}\}
	\end{align*}
	To generalize this result to $z\le x$, we apply \autoref{lm-buchstab} so that
	\begin{align*}
		H(x,z)
		&=A^{-1}(\log x)^kI_k\left({\log x\over\log R},1\right)-\sum_{z\le p<x}h(p)H\left(\frac xp,p\right)+O\{(\log x)^{k-1}\} \\
		&=A^{-1}(\log x)^kI_k\left\{I_k\left({\log x\over\log R},1\right)-\sum_{z\le p<x}\frac kpI_k\left({\log x/p\over\log R},{\log x/p\over\log p}\right)\left(\log x/p\over\log x\right)^k\right\} \\
		&+O\left\{\sum_{z\le p<x}{(\log x)^{k-1}\over p}\right\}+O\left\{\sum_{z\le p<x}{(\log x)^k\over p^2}\right\}+O\{(\log x)^{k-1}\} \\
		&=A^{-1}(\log x)^k\left\{I_k\left({\log x\over\log R},1\right)-k\int_z^x{\mathrm dt\over t\log t}I_k\left({\log x/t\over\log R},{\log x/t\over\log t}\right)\left(\log x/t\over\log x\right)^k\right\} \\
		&+O\{(\log x)^{k-1}\} \\
		&=A^{-1}(\log x)^k\left\{I_k\left({\log x\over\log R},1\right)-k\int_1^{\log x/\log z}I_k\left[(1-\alpha^{-1}){\log x\over\log R},\alpha-1\right]\left(1-\alpha^{-1}\right)^k\mathrm d\alpha\right\} \\
		&+O\{(\log x)^{k-1}\} \\
		&=A^{-1}(\log x)^kI_k\left({\log x\over\log R},{\log x\over\log z}\right)+O\{(\log x)^{k-1}\}.
	\end{align*}
	Plugging these results back into \eqref{eqn-hxz} and \eqref{eqn-qgh}, we obtain
	\begin{align*}
		Q_g\sim{A^2\over(\log R)^{2k}}\cdot A^{-1}(\log R)^kI_k\left({\log R\over\log R},{\log R\over\log z}\right)=A\cdot{I_k(1,u)\over(\log R)^k}.
	\end{align*}
	Combining this result with \eqref{eqn-gp}, we complete the proof.
	\bibliographystyle{plain}
	\bibliography{refs.bib}

\end{document}